\documentclass[12pt, english]{article}

\usepackage{lipsum}
\usepackage{etoolbox}
\makeatletter
\let\@afterindenttrue\@afterindentfalse
\patchcmd\@thm
{\let\thm@indent\indent}{\let\thm@indent\noindent}%
{}{}
\makeatother

\RequirePackage{ifthen}

\newtheorem{thm}{Theorem}[section]
\newtheorem{lemma}[thm]{Lemma}
\newtheorem{corollary}[thm]{Corollary}

\newtheorem{myexamp}{Example}[section]
\newtheorem{definition}[myexamp]{Definition}

\usepackage{amsmath,amssymb,amsfonts}
\usepackage{mathrsfs}
\usepackage{enumitem}

\input{xy}
\xyoption{all}
\RequirePackage[top=4cm,bottom=3.5cm,inner=3.5cm,outer=3.5cm,
headsep=1cm,footskip=1cm]{geometry}

\RequirePackage{xr-hyper}
\RequirePackage{hyperref}

\hypersetup{bookmarks=true}%
\hypersetup{bookmarksopen=true}%

\ifthenelse{\boolean{pdf}}{%
	\hypersetup{colorlinks=true}}{%
	\hypersetup{colorlinks=false}}

\hypersetup{linkcolor=blue}%
\hypersetup{urlcolor=cyan}%
\hypersetup{citecolor=cyan}

\newcommand{\N}{\mathbb{N}}
\newcommand{\K}{\mathbb{K}}
\newcommand{\ideal}[1]{\ensuremath{\mathfrak{#1}}}
\newcommand{\M}[2]{\ensuremath{\mathrm{M}_{#1}(#2)}}
\newcommand{\gen}[2]{\ensuremath{\mathrm{GL}_{#1}(#2)}}
\newcommand{\spl}[2]{\ensuremath{\mathrm{SL}_{#1}(#2)}}
\newcommand{\e}[2]{\ensuremath{\mathrm{E}_{#1}(#2)}}
\newcommand{\re}[3]{\ensuremath{\mathrm{E}_{#1}(#2,#3)}}
\newcommand{\spg}[2]{\ensuremath{\mathrm{Sp}_{#1}(#2)}}
\newcommand{\rspg}[3]{\ensuremath{\mathrm{Sp}_{#1}(#2,#3)}}
\newcommand{\esp}[2]{\ensuremath{\mathrm{ESp}_{#1}(#2)}}
\newcommand{\resp}[3]{\ensuremath{\mathrm{ESp}_{#1}(#2,#3)}}

\begin{document}
	\begin{center}
		\Large{\textbf{Normality theorem for elementary symplectic group with respect to an alternating form}}
	\end{center}
	\begin{center}
		Ruddarraju Amrutha, Pratyusha Chattopadhyay
	\end{center}
		
\medskip

	\begin{center}
		Abstract
	\end{center}

A.A. Suslin proved a normality theorem for an elementary linear group, which says that an elementary linear group of size bigger than or equal to $3$ over a commutative ring with unity is normal in the general linear group of same size. Subsequently, V.I. Kopeiko extended this result of Suslin for a symplectic group defined with respect to the standard skew-symmetric matrix of even size. Here we generalise the result of Kopeiko for a symplectic group defined with respect to any invertible skew-symmetric matrix of even size of Pfaffian one.

\medskip

	\section{Introduction}
	\label{section:1}
	In \cite{Sus}, A.A. Suslin proved a normality theorem for the elementary linear group, which says that for a commutative ring $R$ with $1$, the elementary linear group $\e{n}{R}$ is normal in the general linear group $\gen{n}{R}$, for $n\geq 3$. This normality theorem is important as it was used in proving a $K_1$-analogue of Serre's conjucture. Serre's conjecture says that for a field $\K$, any finitely generated projective module over the polynomial ring $\K[x_1,\ldots,x_r]$ is free. The $K_1$-analogue of this result, proved by Suslin in \cite{Sus}, says that every element in the special linear group $\spl{n}{\K[x_1,\ldots,x_r]}$ can be written as a product of elementary matrices, for $n\geq 3$. Suslin also proved a stronger relative version of the above mentioned normality theorem. 
	
	A similar normality theorem was proved by V.I. Kopeiko in \cite{Kop} for elementary symplectic group. In the symplectic case, the relative version of the normality theorem says that for a commutative ring $R$ with $R=2R$ and an ideal $I$ of $R$, the relative elementary symplectic group $\resp{2n}{R}{I}$ is a normal subgroup of $\spg{2n}{R}$, for $n\geq 2$. Kopeiko used this result to prove that for a field $\K$, any element of $\spg{2n}{\K[x_1,\ldots,x_r]}$ can be written as a product of elementary symplectic matrices over $\K[x_1,\ldots,x_r]$, for $n\geq 2$. In \cite{Kop-Sus}, Suslin and Kopeiko together proved a normality theorem for orthogonal groups which they used to prove that if $\K$ is a field and $A=\K[x_1,\ldots,x_r]$, then any quadratic $A$-space is extended from $\K$.
	
	In this paper, we generalize Kopeiko's normality theorem. The main result (Theorem \ref{lemma:5.2}) says that if $\phi$ is an invertible skew-symmetric matrix of Pfaffian $1$ of size $2n$ over $R$ with $R=2R$, and $I$ is an ideal of $R$, then the relative elementary symplectic group $\resp{\phi}{R}{I}$ with respect to the matrix $\phi$ is normal in $\spg{\phi}{R}$, the symplectic group with respect to $\phi$.   
	\section{Preliminaries}
	\label{section:2}
	Let $R$ be a commutative ring with unity. Let $R^n$ denote the set of columns with $n$ entries in $R$. The ring of matrices of size $n\times n$ with entries in $R$ is denoted by $\M{n}{R}$. The identity matrix of size $n\times n$ is denoted by $I_n$ and $e_{ij}$ denotes the $n\times n$ matrix which has $1$ in the $(i,j)$-th position and $0$ everywhere else. The collection of invertible $n\times n$ matrices with entries in $R$ is denoted by $\gen{n}{R}$. For $\alpha\in\M{m}{R}$ and $\beta\in\M{n}{R}$, the matrix {\tiny{$\begin{pmatrix}\alpha&0\\0&\beta\end{pmatrix}$}}, which is an element of $\M{m+n}{R}$, is denoted by $\alpha\perp\beta$.
	
	\begin{definition}
		The \textit{elementary linear group} $\e{n}{R}$ is the subgroup of $\gen{n}{R}$ generated by elements of the form $E_{ij}(a)=I_n+ae_{ij}$, for $a\in R$. For an ideal $I$ of $R$, the subgroup of $\e{n}{R}$ generated by $E_{ij}(x)$, for $x\in I$ is denoted by $\e{n}{I}$. The \textit{relative elementary group} is the normal closure of $\e{n}{I}$ in $\e{n}{R}$. It is denoted by $\re{n}{R}{I}$.  
	\end{definition}
	
	\begin{definition}
		Given an invertible skew-symmetric matrix $\phi$ of size $2n\times 2n$, 
		\begin{equation*} \spg{\phi}{R}=\{\alpha\in\gen{2n}{R}\;\big|\;\alpha^t\phi\alpha=\phi\}.
		\end{equation*}
		
		Let $\psi_n=\sum_{i=1}^{n}(e_{2i-1,2i}-e_{2i,2i-1})$ denote the standard skew-symmetric matrix. The \textit{symplectic group}, denoted by $\spg{2n}{R}$, is a subgroup of $\gen{2n}{R}$ given by 
		\begin{equation*} \spg{2n}{R}=\{\alpha\in\gen{2n}{R}\;\big|\;\alpha^t\psi_n\alpha=\psi_n\}.
		\end{equation*}
		
	\end{definition}
	
	\begin{definition}
		Let $I$ be an ideal of $R$. Then we have the quotient map $f:R\rightarrow R/I$. Using $f$, we can define a ring homomorphism $\tilde{f}:\spg{\phi}{R}\rightarrow\spg{\phi}{R/I}$ given by $\tilde{f}(a_{ij})=(f(a_{ij}))$. We denote the kernel of this map by $\rspg{\phi}{R}{I}$.
	\end{definition}
	
	\begin{definition}
		The \textit{elementary symplectic group}, denoted by $\esp{2n}{R}$, is a subgroup of $\spg{2n}{R}$ generated by elements of the form $se_{ij}(a)$, for $a\in R$ and
		\begin{equation*}
			se_{ij}(a)=\begin{cases}
				I_n+ae_{ij},&\text{ if } i=\sigma(j)\\
				I_n+ae_{ij}-(-1)^{i+j}ae_{\sigma(j)\sigma(i)},&\text{ if } i\neq\sigma(j)\\	
			\end{cases}
		\end{equation*} where $\sigma$ denotes the permutation of $\{1,\ldots,2n\}$ given by $\sigma(2i-1)=2i$ and $\sigma(2i)=2i-1$. We denote by $\esp{2n}{I}$ a subgroup of $\esp{2n}{R}$ generated by elements of the form $se_{ij}(x)$, for $x\in I$. The \textit{relative elementary group}, denoted by $\resp{2n}{R}{I}$, is the normal closure of $\esp{2n}{I}$ in $\esp{2n}{R}$.
	\end{definition}

	\noindent
	\textbf{Remark:} For a ring $R$ with $R=2R$ and an ideal $I$ of $R$, the relative elementary group $\resp{2n}{R}{I}$ is the smallest normal subgroup of $\esp{2n}{R}$ containing $se_{21}(x)$, for all $x\in I$. A similar result for the elementary linear group was proved by W. van der Kallen in \cite{van} (Lemma 2.2). The result for symplectic group can be proved using a similar argument.
	
	\begin{definition}
		Let $\phi$ be an invertible skew-symmetric matrix of size $2n$ of the form $\begin{pmatrix} 0 & -c^t\\ c & \nu \end{pmatrix}$, and $\phi^{-1}$ be of the form $\begin{pmatrix} 0 & d^t\\ -d & \mu \end{pmatrix}$, where $c,d\in R^{2n-1}$.
		
		Given $v\in R^{2n-1}$, consider the matrices $\alpha$ and $\beta$ defined as 
		\begin{equation*}
			\begin{aligned}
				\alpha&:=\alpha_\phi(v)&:=I_{2n-1}+dv^t\nu\\
				\beta&:=\beta_\phi(v)&:=I_{2n-1}+\mu vc^t.
			\end{aligned}
		\end{equation*}
		L.N. Vaserstein constructed these matrices in Lemma 5.9, \cite{Vas}. Note that $\alpha$ and $\beta$ depend on $\phi$ and $v$. Also, $\alpha,\beta\in\e{2n-1}{R}$ by Lemma 2.2, \cite{Vas}. Using these matrices, Vaserstein constucted the following matrices in \cite{Vas}:
		\begin{equation*}
			\begin{aligned}
				C_\phi(v)&:=\begin{pmatrix} 1&0\\ v&\alpha\end{pmatrix}\\
				R_\phi(v)&:=\begin{pmatrix} 1&v^t\\ 0&\beta\end{pmatrix}.
			\end{aligned}
		\end{equation*}
		
		The \textit{elementary symplectic group} $\esp{\phi}{R}$ \textit{with respect to the invertible skew-symmetric matrix} $\phi$ is a subgroup of $\spg{\phi}{R}$ generated by $C_\phi(v)$ and $R_\phi(v)$, for $v\in R^{2n-1}$. We denote by $\esp{\phi}{I}$ a subgroup of $\esp{\phi}{R}$ generated as a group by the elements $C_\phi(v)$ and $R_\phi(v)$, for $v\in I^{2n-1}(\subseteq R^{2n-1})$. The \textit{relative elementary symplectic group} $\resp{\phi}{R}{I}$ is the normal closure of $\esp{\phi}{I}$ in $\esp{\phi}{R}$.
	\end{definition}
	
		%
		%
		%
	
	\section{Results about elementary symplectic group}
	\label{section:3}
	In this section, we will recall a few results related to elementary symplectic groups. We will also establish normality theorem for the elementary symplectic group over a local ring (see Lemma \ref{lemma:3.8}) and over a polynomial ring in one variable over a local ring (see Lemma \ref{lemma:3.10}). Finally, we will establish Corollary \ref{lemma:3.13} which plays a crucial role in the proof of the main theorem (Theorem \ref{lemma:5.2}).
	 
	\begin{lemma}(Corollary 1.11, \cite{Kop})
		\label{lemma:3.1}
		Let $R$ be a ring and $I$ be an ideal of $R$. Let $n\geq 2$. Then, $\resp{2n}{R}{I}$ is a normal subgroup of $\spg{2n}{R}$.
	\end{lemma}
	
	\begin{lemma}(Remark 4.3, \cite{JOA})
		\label{lemma:3.2}
		Let $(R,\ideal{m})$ be a local ring and $\phi$ be a skew-symmetric matrix of Pfaffian $1$ of size $2n$ over $R$. Then $\phi=\epsilon^t\psi_n\epsilon$, for some $\epsilon\in\e{2n}{R}$. 
	\end{lemma}
	
	\begin{lemma}(\cite{RaoSwan})
		\label{lemma:3.3}
		For $n\geq 2$ and $\epsilon\in\e{2n}{R}$, we have an $\epsilon_0\in\e{2n-1}{R}$ such that $\epsilon^t\psi_n\epsilon= (1\perp\epsilon_0)^t\psi_n(1\perp\epsilon_0)$.
	\end{lemma}
	
	\begin{lemma}(Lemma 3.6, 3.7, \cite{JOA})
		\label{lemma:3.4}
		Let $\phi$ and $\phi^\ast$ be two invertible skew-symmetric matrices such that $\phi=(1\perp\epsilon)^t\phi^\ast(1\perp\epsilon)$ for some $\epsilon\in\e{2n-1}{R}$. Then, we have
		\begin{equation*}
			\begin{aligned}
				\spg{\phi}{R}&=(1\perp\epsilon)^{-1}\spg{\phi^\ast}{R}(1\perp\epsilon),\\
				\esp{\phi}{R}&=(1\perp\epsilon)^{-1}\esp{\phi^\ast}{R}(1\perp\epsilon).
			\end{aligned}
		\end{equation*}
	\end{lemma}
	
	\begin{lemma}(Lemma 3.8, \cite{JOA})
		\label{lemma:3.5}
		Let $\phi$ and $\phi^\ast$ be two invertible skew-symmetric matrices such that $\phi=(1\perp\epsilon)^t\phi^\ast(1\perp\epsilon)$, for some $\epsilon\in\re{2n-1}{R}{I}$. Then,
		\begin{equation*}
			\resp{\phi}{R}{I}=(1\perp\epsilon)^{-1}\resp{\phi^\ast}{R}{I}(1\perp\epsilon).
		\end{equation*}
	\end{lemma}
	
	The following lemma shows that the elementary symplectic group $\esp{\phi}{R}$ with respect to a skew-symmetric matrix $\phi$ can be considered as a generalisation of the elementary symplectic group $\esp{2n}{R}$. We include the proof for completeness.
	\begin{lemma}(Lemma 3.5, \cite{JOA})
		\label{lemma:3.6}
		Let $R$ be a ring with $R=2R$ and $n\geq 2$. Then,  $\esp{\psi_n}{R}=\esp{2n}{R}$.
	\end{lemma}
	\textit{Proof:} $\esp{\psi_n}{R}\subseteq\esp{2n}{R}$ as for $v=(a_1,\cdots,a_{2n-1})^t\in R^{2n-1}$, we have 
	\begin{equation*}
		C_{\psi_n}(v)=\prod_{i=2}^{2n}se_{i1}(a_{i-1}) \text{ and }
		R_{\psi_n}(v)=\prod_{i=2}^{2n}se_{1i}(a_{i-1})
	\end{equation*}
	For integers $i,j$ with $i\neq j,\sigma(j)$ and for $a,b\in R$, we have 
	\begin{equation*}
		\begin{aligned}
			[se_{i \sigma(i)}(a),se_{\sigma(i) j}(b)]&=se_{ij}(ab)se_{\sigma(j)j}((-1)^{i+j}ab^2),\\
			[se_{ik}(a),se_{kj}(b)]&=se_{ij}(ab), \text{ if  } k\neq\sigma(i),\sigma(j),\\
			[se_{ik}(a),se_{k\sigma(i)}(b)]&=se_{i\sigma(i)}(2ab), \text{ if } k\neq i,\sigma(i).
		\end{aligned}
	\end{equation*}
	Using these identities, $se_{ij}(a)$, for $i,j\neq 1$, can be written as a product of elements of the form $se_{1i}(x)$ and $se_{j1}(y)$, for $x,y\in R$. Also, $se_{1i}(a), se_{j1}(b)\in\esp{\psi_n}{R}$. So, $\esp{2n}{R}\subseteq\esp{\psi_n}{R}$. \\
	
	\noindent A relative version of the above lemma with respect to an ideal is as follows.
	
	\begin{lemma}(Lemma 3.5, \cite{JOA})
		\label{lemma:3.7}
		Let $R$ be a ring with $R=2R$ and $n\geq 2$. For an ideal $I$ of $R$, we have $\resp{\psi_n}{R}{I}=\resp{2n}{R}{I}$. 
	\end{lemma}
	\textit{Proof:} $\resp{\psi_n}{R}{I}$ is generated by elements of the form $\gamma\delta(v)\gamma^{-1}$, where $\gamma\in\esp{\psi_n}{R}$ and $\delta(v)$  is either $C_{\psi_n}(v)$ or $R_{\psi_n}(v)$, with $v\in I^{2n-1}$. Since $C_{\psi_n}((a_1,\ldots,a_n)^t)=\prod_{i=2}^{2n}se_{i1}(a_{i-1})$ and $R_{\psi_n}((a_1,\ldots,a_n)^t)=\prod_{i=2}^{2n}se_{1i}(a_{i-1})$, we have $\delta(v)\in\esp{2n}{I}$. By Lemma \ref{lemma:3.6}, $\esp{\psi_n}{R}=\esp{2n}{R}$ and hence $\gamma\in\esp{2n}{R}$. By definition of $\resp{2n}{R}{I}$, we have $\gamma\delta(v)\gamma^{-1}\in\resp{2n}{R}{I}$. Hence $\resp{\psi_n}{R}{I}\subseteq\resp{2n}{R}{I}$.
	
	$\resp{2n}{R}{I}$ is the smallest normal subgroup of $\esp{2n}{R}$ containing $se_{21}(x)$, for all $x\in I$ (see remark in section \ref{section:1}). For $x\in I$, we have $se_{21}(x)=C_{\psi_n}((x,0,\ldots,0)^t)$. So, for $\gamma\in\esp{2n}{R}=\esp{\psi_n}{R}$, we have $\gamma se_{21}(x)\gamma^{-1}\in\resp{\psi_n}{R}{I}$ and hence we have the other way inclusion. Therefore, $\resp{\psi_n}{R}{I}=\resp{2n}{R}{I}$.
	
	\begin{lemma}
		\label{lemma:3.8}
		Let $(R,\ideal{m})$ be a local ring. Suppose that $R=2R$. Let $\phi$ be a skew-symmetric matrix of Pfaffian $1$ of size $2n$ over $R$ with $n\geq 2$. Then, $\esp{\phi}{R}$ is a normal subgroup of $\spg{\phi}{R}$.
	\end{lemma}
	\textit{Proof:} Let $\phi$ be a skew-symmetric matrix of Pfaffian $1$ of size $2n$ over $R$. By Lemma \ref{lemma:3.2} and Lemma \ref{lemma:3.3}, we can write $\phi=(1\perp\epsilon)^t\psi_n(1\perp\epsilon)$, for some $\epsilon\in\e{2n-1}{R}$. Hence, by Lemma \ref{lemma:3.4} and Lemma \ref{lemma:3.6}, we have $\spg{\phi}{R}=(1\perp\epsilon)^{-1}\spg{2n}{R}(1\perp\epsilon)$ and  $\esp{\phi}{R}=(1\perp\epsilon)^{-1}\esp{2n}{R}(1\perp\epsilon)$.
	
	Let $\gamma\in\spg{\phi}{R}$ and $\delta\in\esp{\phi}{R}$. Then $\delta=(1\perp\epsilon)^{-1}\delta_1(1\perp\epsilon)$ and $\gamma=(1\perp\epsilon)^{-1}\gamma_1(1\perp\epsilon)$ for some $\delta_1\in\esp{2n}{R}$ and $\gamma_1\in\spg{2n}{R}$. Hence $\gamma\delta\gamma^{-1}=(1\perp\epsilon)^{-1}\gamma_1(1\perp\epsilon)(1\perp\epsilon)^{-1}\delta_1(1\perp\epsilon)(1\perp\epsilon)^{-1}\gamma_1^{-1}(1\perp\epsilon)= (1\perp\epsilon)^{-1}\gamma_1\delta_1\gamma_1^{-1}(1\perp\epsilon)$. By Lemma \ref{lemma:3.1}, $\gamma_1\delta_1\gamma_1^{-1}\in\esp{2n}{R}$ and hence $\gamma\delta\gamma^{-1}\in\esp{\phi}{R}$. Therefore, $\esp{\phi}{R}$ is a normal subgroup of $\spg{\phi}{R}$.\\
	
	\noindent A relative version of this result (see Lemma \ref{lemma:3.9}) can be proved in a similar manner.
	
	\begin{lemma}
		\label{lemma:3.9}
		Let $(R,\ideal{m})$ be a local ring. Suppose that $R=2R$. Let $\phi$ be a skew-symmetric matrix of Pfaffian $1$ of size $2n$ over $R$ with $n\geq 2$. For an ideal $I$ of $R$, we have $\resp{\phi}{R}{I}$ is a normal subgroup of $\spg{\phi}{R}$.
	\end{lemma} 
	
	\noindent
	\textbf{Notation:} Let $\phi$ be an invertible skew-symmetric matrix of size $2n$ over $R$. 
	\begin{equation*}
		\spg{\phi\otimes R[X]}{R[X]}:=\{\alpha\in\gen{2n}{R[X]}\;\big|\; \alpha^t\phi\alpha=\phi\}.
	\end{equation*}
	By $\esp{\phi\otimes R[X]}{R[X]}$, we mean the elementary symplectic group generated by $C_\phi(v)$ and $R_\phi(v)$, where $v\in R[X]^{2n-1}$.\\
	
	 Let $R$ and $S$ be two commutative rings with unity and let $f:R\rightarrow S$ be a ring homomorphism. If $\phi=(a_{ij})$, then define $f(\phi)=(f(a_{ij}))$. For a maximal ideal $\ideal{m}$ of $R$, the group $\esp{\phi\otimes R_\ideal{m}[X]}{R_\ideal{m}[X]}$ is the elementary symplectic group with respect to the matrix $f(\phi)$, where $f:R\rightarrow R_\ideal{m}[X]$ is the map $a\mapsto \frac{a}{1}$.   
	
	\begin{lemma}
		\label{lemma:3.10}
		Let $\phi$ be an invertible skew-symmetric matrix of Pfaffian 1 of size $2n$ over $R$ with $R=2R$. Let $\ideal{m}$ be a maximal ideal of $R$. Then, $\esp{\phi\otimes R_\ideal{m}[X]}{R_\ideal{m}[X]}$ is a normal subgroup of $\spg{\phi\otimes R_\ideal{m}[X]}{R_\ideal{m}[X]}$. 
	\end{lemma}
	\textit{Proof:} Consider the local ring $(R_\ideal{m},\ideal{m})$. Then $\phi$ is a skew-symmetric matrix of Pfaffian 1 over $R_\ideal{m}$. By Lemma \ref{lemma:3.2} and Lemma \ref{lemma:3.3}, there exists $\epsilon\in\e{2n-1}{R_\ideal{m}[X]}$ such that $\phi=(1\perp\epsilon)^t\psi_n(1\perp\epsilon)$. Hence by Lemma \ref{lemma:3.4} and Lemma \ref{lemma:3.6}, we can write $\spg{\phi}{R_\ideal{m}[X]}=(1\perp\epsilon)^{-1}\spg{2n}{R_\ideal{m}[X]}(1\perp\epsilon)$ and $\esp{\phi}{R_\ideal{m}[X]}=(1\perp\epsilon)^{-1}\esp{2n}{R_\ideal{m}[X]}(1\perp\epsilon)$. 
	
	Since $\esp{2n}{R_\ideal{m}[X]}$ is normal in $\spg{2n}{R_\ideal{m}[X]}$ by Lemma \ref{lemma:3.1}, it follows that $\esp{\phi}{R_\ideal{m}[X]}$ is normal in $\spg{\phi}{R_\ideal{m}[X]}$.
	
	
	\begin{lemma}
		\label{lemma:3.11}
		Let $\phi$ be an invertible skew-symmetric matrix of size $2n$ over $R$ with $R=2R$. For $v,w\in R^{2n-1}$, we have the following splitting property for the generators of $\esp{\phi}{R}$:  
		\begin{equation*}
			\begin{aligned}
				C_\phi(v+w)&=C_\phi(v/2)C_\phi(w)C_\phi(v/2),\\
				R_\phi(v+w)&=R_\phi(v/2)R_\phi(w)R_\phi(v/2).
			\end{aligned}
		\end{equation*}
	\end{lemma}
	
	\noindent \textit{Proof:} Suppose $\phi$ is of the form $\begin{pmatrix} 0 & -c^t\\ c & \nu \end{pmatrix}$ and $\phi^{-1}$ is of the form $\begin{pmatrix} 0 & d^t\\ -d & \mu \end{pmatrix}$. Then,
	\begin{equation*}
		\begin{aligned}
			\alpha_\phi(v)\alpha_\phi(w)&=(I_{2n-1}+dv^t\nu)(I_{2n-1}+dw^t\nu)\\ 
			&=I_{2n-1}+dw^t\nu+dv^t\nu\\
			&=I_{2n-1}+d(v+w)^t\nu= \alpha_\phi(v+w).
		\end{aligned}
	\end{equation*} 
	Similarly, $\beta_\phi(v+w)=\beta_\phi(v)\beta_\phi(w)$. Now,  
	\begin{eqnarray*}
			&&C_\phi(v/2)C_\phi(w)C_\phi(v/2)\\
			&=&\begin{pmatrix}1&0\\v/2&\alpha_\phi(v/2)\end{pmatrix}
			\begin{pmatrix}1&0\\w&\alpha_\phi(w)\end{pmatrix}
			\begin{pmatrix}1&0\\v/2&\alpha_\phi(v/2)\end{pmatrix}\\
			&=& \begin{pmatrix}1&0\\(v/2+\alpha_\phi(v/2)w+\alpha_\phi(v/2)\alpha_\phi(w)v/2)&(\alpha_\phi(v/2)\alpha_\phi(w)\alpha_\phi(v/2))\end{pmatrix}\\
			&=&\begin{pmatrix}1&0\\v+w&\alpha_\phi(v+w)\end{pmatrix}= C_\phi(v+w).
	\end{eqnarray*}
	Similarly, $R_\phi(v+w)=R_\phi(v/2)R_\phi(w)R_\phi(v/2)$. 
	
	\begin{lemma}
		\label{lemma:2.14}
		Let $G$ be a group and $a_i,b_i\in G$ for $i-1,\ldots,n$. Then 
		\begin{equation*}
			\prod_{i=1}^na_ib_i=\bigg(\prod_{i=1}^n\big(\prod_{j=1}^ia_i\big)b_i\big(\prod_{j=1}^ia_i\big)^{-1}\bigg)\prod_{i=1}^na_i.
		\end{equation*}
	\end{lemma}
	\begin{lemma}
		\label{lemma:3.12}
		Let $R$ be a ring with $R=2R$ and $\phi$ be an invertible skew-symmetric matrix of size $2n$. For an ideal $I$ of $R$, we have $\resp{\phi}{R}{I}=\esp{\phi}{R}\cap\rspg{\phi}{R}{I}$. 
	\end{lemma}
	\textit{Proof:} Elements of $\resp{\phi}{R}{I}$ are generated by elements of the form $\gamma(v)\delta(w)\gamma(v)^{-1}$, where $\gamma(v)$ denotes $C_\phi(v)$ or $R_\phi(v)$ and $\delta(w)$ denotes $C_\phi(w)$ or $R_\phi(w)$, where $v\in R^{2n-1}$ and $w\in I^{2n-1}$. Also, $\gamma(v)\delta(w)\gamma(v)^{-1}=I_{2n} (mod\text{ }I)$. So, $\resp{\phi}{R}{I}\subseteq\esp{\phi}{R}\cap\rspg{\phi}{R}{I}$.  
	
	Let $\gamma=\prod_{i=1}^t\gamma_i(v_i)\in\esp{\phi}{R}\cap\rspg{\phi}{R}{I}$ where $\gamma_i(v_i)$ denotes $C_\phi(v_i)$ or $R_\phi(v_i)$ with $v_i\in R^{2n-1}$. Then we have $\prod_{i=1}^t\gamma_i(v_i) (mod \text{ }I)=I_{2n}$ which implies that for each $i$, there exist $u_i\in R^{2n-1}$ and $w_i\in I^{2n-1}$ such that $\gamma=\prod_{i=1}^t\gamma_i(u_i+w_i)$ with $\gamma_i(u_i)=I_{2n}$. Now, 
	\begin{equation*}
		\begin{aligned}
			\gamma&=\prod_{i=1}^t\gamma_i(u_i+w_i)\\
			&=\prod_{i=1}^t(\gamma_i(u_i/2)\gamma_i(w_i)\gamma_i(u_i/2)) \text{ [by Lemma \ref{lemma:3.11}}]\\
			&=\bigg(\prod_{i=1}^{t-1}\eta_i\zeta_i\bigg)\eta_n, \text{ for some } \eta_i\in \resp{\phi}{R}{I} \text{ and } \zeta_i\in\esp{\phi}{R}.\\
			&(\text{Take } \eta_i=\gamma_i(\pm u_i/2)\gamma_i(w_i)\gamma_i(\mp u_i/2)) \text{ and } \zeta_i=\gamma_i(u_i)\gamma_{i+1}(u_{i+1})).\\
			&= \eta_1 \bigg(\prod_{i=1}^t \bigg(\prod_{j=1}^i\zeta_j\bigg)\eta_{i-1} \bigg(\prod_{j=1}^i\zeta_j\bigg)^{-1}\bigg) \bigg(\prod_{i=1}^{t-1}\zeta_i\bigg)\eta_n \text{ [by Lemma \ref{lemma:2.14}]}\\
			&= \eta_1 \bigg(\prod_{i=1}^{t-1} \bigg(\prod_{j=1}^i\zeta_j\bigg)\eta_{i-1} \bigg(\prod_{j=1}^i\zeta_j\bigg)^{-1}\bigg) \eta_n \in\resp{\phi}{R}{I},
		\end{aligned}	
	\end{equation*}
	and hence $\esp{\phi}{R}\cap\rspg{\phi}{R}{I}\subseteq\resp{\phi}{R}{I}$. Therefore, we have the required equality. 
	
	\begin{corollary}
		\label{lemma:3.13}
		$\resp{\phi}{R[X]}{(X)}=\esp{\phi}{R[X]}\cap\rspg{\phi}{R[X]}{(X)}$.
	\end{corollary}
	\textit{Proof:}
	Follows from Lemma \ref{lemma:3.12}.

	\section{Local-Global principle for $\esp{\phi\otimes R[X]}{R[X]}$}	
	\label{section:4}
	We will state the Local-Global principle for $\esp{\phi\otimes R[X]}{R[X]}$ which was proved in \cite{JOA}. Using this, we will prove a graded case of the Local-Global principle for $\esp{\phi\otimes R[X_1,\ldots,X_t]}{R[X_1,\ldots,X_t]}$.  
	
	\begin{lemma}{(Local-Global principle)}(Theorem 4.6, \cite{JOA})
		\label{lemma:4.1} 
		Let $\phi$ be a skew-symmetric matrix of Pfaffian $1$ of size $2n$ over $R$ with $n\geq 2$. Let $\theta(X)\in\spg{\phi\otimes R[X]}{R[X]}$, with $\theta(0)=I_{2n}$. If $\theta(X)_\ideal{m}\in\esp{\phi\otimes R_\ideal{m}[X]}{R_\ideal{m}[X]}$, for all maximal ideals $\ideal{m}$ of $R$, then $\theta(X)\in\esp{\phi\otimes R[X]}{R[X]}$. 
	\end{lemma}
	
	Next, we will prove the graded case of the Local-Global principle. To prove this result, we borrow ideas from Lemma 3.1 of \cite{Rao}.
	\begin{lemma}{(Graded case of Local-Global principle)}
		\label{lemma:4.2}
		Let $\phi$ be a skew-symmetric matrix of Pfaffian $1$ of size $2n$ over $R$ with $n\geq 2$. Let 
		\begin{equation*}
			\theta(X_1,\ldots,X_t)\in\spg{\phi\otimes R[X_1,\ldots,X_t]}{R[X_1,\ldots,X_t]},
		\end{equation*} with $\theta(0,\ldots,0)=I_{2n}$. If $\theta(X_1,\ldots,X_t)_\ideal{m}\in\esp{\phi\otimes R_\ideal{m}[X_1,\ldots,X_t]}{R_\ideal{m}[X_1,\ldots,X_t]}$, for all maximal ideals $\ideal{m}$ of $R$, then $\theta(X_1,\ldots,X_t)\in\esp{\phi\otimes R[X_1,\ldots,X_t]}{R[X_1,\ldots,X_t]}$. 
	\end{lemma}
	\textit{Proof:}
	Let $S$ denote the polynomial ring $R[X_1,\ldots,X_t]$. Let $S_0=R$ and for $i\in\N$, let $S_i$ be the collection of homogeneous polynomials in $X_1,\ldots,X_t$ of degree $i$. Then, $S$ can be considered as a graded ring with the grading $S=S_0+S_1+S_2+\cdots$. Any element of $S$ can be written uniquely as $a_0+a_1+a_2+\cdots$, where $a_i\in S_i$. Define a ring homomorphism $f:S\rightarrow S[T]$ as $f(a_0+a_1+a_2+\cdots)=a_0+a_1T+a_2T^2+\cdots$. Define $\tilde{\theta}(T)$ to be the matrix obtained by taking the image of $f$ on each entry of $\theta(X_1,\ldots,X_t)$. That is, if $\theta(X_1,\ldots,X_t)=(a_{ij})$, then  $\tilde{\theta}(T)=(f(a_{ij}))$. Then, $\tilde{\theta}(T)\in\spg{\phi\otimes S[T]}{S[T]}$ as $\theta(X_1,\ldots,X_t)\in\spg{\phi\otimes R[X_1,\ldots,X_t]}{R[X_1,\ldots,X_t]}$. Also, $\tilde{\theta}(0)=\theta(0,\ldots,0)=I_{2n}$ and $\tilde{\theta}(1)=\theta(X_1,\ldots,X_t)$
	
	Let $\ideal{m}$ be a maximal ideal of $R$. By the hypothesis, 
	\begin{equation*} 
		\theta(X_1,\ldots,X_t)_\ideal{m}\in\esp{\phi\otimes R_\ideal{m}[X_1,\ldots,X_t]}{R_\ideal{m}[X_1,\ldots,X_t]}.
	\end{equation*} 
	There exists $s_\ideal{m}\in R\setminus\ideal{m}$ such that $\theta(X_1,\ldots,X_t)_{s_\ideal{m}}\in\esp{\phi\otimes R_{s_\ideal{m}}[X_1,\ldots,X_t]}{R_{s_\ideal{m}}[X_1,\ldots,X_t]}$. This implies that $\tilde{\theta}(T)_{s_\ideal{m}}\in\esp{\phi\otimes S_{s_\ideal{m}}[T]}{S_{s_\ideal{m}}[T]}$. Now, for every maximal ideal $\ideal{m}$ of $R$, we have $s_\ideal{m}\in R\setminus\ideal{m}$. So, the set $\{s_\ideal{m}\,|\,\ideal{m}\text{ is a maximal ideal of }R\}$  generates $R$ and hence generates $S$.
	
	Let $\ideal{m}'$ be a maximal ideal of $S$. Since $\{s_\ideal{m}\,|\,\ideal{m}\text{ is a maximal ideal of }R\}$ generates $S$ and $\ideal{m}'$ is a proper ideal of $S$, there exists a maximal ideal $\ideal{m}$ of $R$ such that $s_\ideal{m}\notin\ideal{m}'$. By choice of $s_\ideal{m}$,  $\tilde{\theta}(T)_{s_\ideal{m}}\in\esp{\phi\otimes S_{s_\ideal{m}}[T]}{S_{s_\ideal{m}}[T]}$. So, $\tilde{\theta}(T)_{\ideal{m'}}\in\esp{\phi\otimes S_{\ideal{m'}}[T]}{S_{\ideal{m'}}[T]}$. By Lemma \ref{lemma:4.1}, we have $\tilde{\theta}(T)\in\esp{\phi\otimes S[T]}{S[T]}$. Substituting $T=1$, we have $\theta(X_1,\ldots,X_t)\in\esp{\phi\otimes R[X_1,\ldots,X_t]}{R[X_1,\ldots,X_t]}$.

	\section{Normality of $\resp{\phi}{R}{I}$}
	\label{section:5} 
	In this section, we prove the main theorem, namely normality of $\resp{\phi}{R}{I}$ (see Theorem \ref{lemma:5.2}). This can be considered as a generalization of famous normality theorem for the elementary symplectic group proved by Kopeiko (see Corollary 1.11 of \cite{Kop}). In Theorem \ref{lemma:5.2}, if we substitute $\phi$ by the standard skew-symmetric matrix $\psi_n$, we get Kopeiko's normality theorem for th elementary symplectic group. 
	\begin{lemma}
		\label{lemma:5.1}
		Let $R$ be a ring with $R=2R$. Let $\phi$ be a skew-symmetric matrix of Pfaffian $1$ of size $2n$ over $R$ with $n\geq 2$. Then, $\esp{\phi}{R}$ is a normal subgroup of $\spg{\phi}{R}$.
	\end{lemma}
	\textit{Proof:}
	Let $\gamma\in\spg{\phi}{R}$ and $\delta\in\esp{\phi}{R}$. Since $\esp{\phi}{R}$ is generated by elements of the form $C_\phi(v)$ and $R_\phi(v)$ for $v\in R$, we can write $\delta=\prod D_i(v_i)$, where $D_i(v_i)$ denotes $C_\phi(v_i)$ or $R_\phi(v_i)$ with $v_i\in R^{2n-1}$.
	
	Define $\delta(X)=\prod D_i(v_iX)$. Then, $\delta(X)\in\esp{\phi\otimes R[X]}{R[X]}$ with $\delta(0)=I_{2n}$ and $\delta(1)=\delta$. Define $\theta(X)=\gamma\delta(X)\gamma^{-1}$. Note that $\theta(X)\in\spg{\phi\otimes R[X]}{R[X]}$ with $\theta(0)=I_{2n}$.
	
	Let $\ideal{m}$ be a maximal ideal of $R$. Then, by Lemma \ref{lemma:3.10}, $\esp{\phi\otimes R_\ideal{m}[X]}{R_\ideal{m}[X]}$ is a normal subgroup of $\spg{\phi\otimes R_\ideal{m}[X]}{R_\ideal{m}[X]}$ and hence $\theta(X)_\ideal{m}=\gamma_\ideal{m}\delta(X)_\ideal{m}\gamma_\ideal{m}^{-1}\in \esp{\phi\otimes R_\ideal{m}[X]}{R_\ideal{m}[X]}$. Since this is true for all maximal ideals $\ideal{m}$ of $R$, by Lemma \ref{lemma:4.1}, we have $\theta(X)\in \esp{\phi\otimes R[X]}{R[X]}$, i.e., $\gamma\delta(X)\gamma^{-1}\in \esp{\phi\otimes R[X]}{R[X]}$. Substituting $X=1$, we get $\gamma\delta\gamma^{-1}\in \esp{\phi}{R}$. 
	
	Therefore, we have $\gamma\delta\gamma^{-1}\in \esp{\phi}{R}$, for all $\gamma\in\spg{\phi}{R}$ and $\delta\in\esp{\phi}{R}$ and hence $\esp{\phi}{R}$ is a normal subgroup of $\spg{\phi}{R}$.

	\begin{thm}
		\label{lemma:5.2}
		Let $\phi$ be a skew-symmetric matrix of Pfaffian $1$ of size $2n$ over $R$ with $n\geq 2$. Assume that $R=2R$. Let $I$ be an ideal of $R$. Then, $\resp{\phi}{R}{I}$ is a normal subgroup of $\spg{\phi}{R}$.
	\end{thm}
	\textit{Proof:}
	Let $\mu\in\spg{\phi}{R}$ and $\lambda\in\resp{\phi}{R}{I}$. Since $\resp{\phi}{R}{I}$ is generated by elements of the form $\gamma(v)\delta(w)(\gamma(v))^{-1}$, where $\gamma(v)$ denotes $C_\phi(v)$ or $R_\phi(v)$ and $\delta(w)$ denotes $C_\phi(w)$ or $R_\phi(w)$ with $v\in R^{2n-1}$ and $w\in I^{2n-1}$, we can write $\lambda=\prod_{i=1}^t\gamma(v_i)\delta(w_i)(\gamma(v_i))^{-1}$, for some $v_i\in R^{2n-1}$ and $w_i\in I^{2n-1}$.
	
	Let $\theta(X_1,\ldots,X_t)=\mu(\prod_{i=1}^t\gamma(v_i)\delta(X_i)(\gamma(v_i))^{-1}) \mu^{-1}$. Then, $\theta(X_1,\ldots,X_t)\in\spg{\phi\otimes R[X_1,\ldots,X_t]}{R[X_1,\ldots,X_t]}$ with $\theta(0,\ldots,0)=I_{2n}$ and $\theta(w_1,\ldots,w_t)=\mu\lambda\mu^{-1}$.
	
	Let $\ideal{m}$ be a maximal ideal of $R$. Then 
	\begin{equation*}
		\theta(X_1,\ldots,X_t)_\ideal{m}= \mu_\ideal{m}\bigg(\prod_{i=1}^t\gamma(v_i)\delta(X_i)(\gamma(v_i))^{-1}\bigg)_\ideal{m} \mu_\ideal{m}^{-1}.	
	\end{equation*}
	Here, $\mu_\ideal{m}\in\spg{\phi\otimes R_\ideal{m}[X_1,\ldots,X_t]}{R_\ideal{m}[X_1,\ldots,X_t]}$ and $(\prod_{i=1}^t\gamma(v_i)\delta(X_i)(\gamma(v_i))^{-1})_\ideal{m}\in\esp{\phi\otimes R_\ideal{m}[X_1,\ldots,X_t]}{R_\ideal{m}[X_1,\ldots,X_t]}$. By Lemma \ref{lemma:5.1}, $\esp{\phi\otimes R_\ideal{m}[X_1,\ldots,X_t]}{R_\ideal{m}[X_1,\ldots,X_t]}$ is a normal subgroup of $\spg{\phi\otimes R_\ideal{m}[X_1,\ldots,X_t]}{R_\ideal{m}[X_1,\ldots,X_t]}$. Therefore, we have 
	\begin{equation*}
		\theta(X_1,\ldots,X_t)_\ideal{m}\in\esp{\phi\otimes R_\ideal{m}[X_1,\ldots,X_t]}{R_\ideal{m}[X_1,\ldots,X_t]}.	
	\end{equation*} 
	This is true for all maximal ideals $\ideal{m}$ of $R$. Hence, by Lemma \ref{lemma:4.2}, we have $\theta(X_1,\ldots,X_t)\in\esp{\phi\otimes R[X_1,\ldots,X_t]}{R[X_1,\ldots,X_t]}$.
	
	Also, $\theta(X_1,\ldots,X_t)\in\rspg{\phi\otimes R[X_1,\ldots,X_t]}{R[X_1,\ldots,X_t]}{(X_1,\ldots,X_t)}$. Hence by Lemma \ref{lemma:3.13}, we get $\theta(X_1,\ldots,X_t)\in\resp{\phi\otimes R[X_1,\ldots,X_t]}{R[X_1,\ldots,X_t]}{(X_1,\ldots,X_t)}$. Substituting $X_i=w_i$ for $i=1,\ldots,t$, we have $\theta(w_1,\ldots,w_t)\in\resp{\phi}{R}{I}$, that is, $\mu\lambda\mu^{-1}\in\resp{\phi}{R}{I}$.
	
	Therefore, $\mu\lambda\mu^{-1}\in\resp{\phi}{R}{I}$ for all $\mu\in\spg{\phi}{R}$ and $\lambda\in\resp{\phi}{R}{I}$ and hence $\resp{\phi}{R}{I}$ is a normal subgroup of $\spg{\phi}{R}$.

\end{document}